\documentclass[11pt]{article}
\usepackage{amsfonts}
\usepackage[dvips]{graphicx}
\newcommand{\halmos}{\rule{1ex}{1.4ex}}
\newcommand{\proofbox}{\hspace*{\fill}\mbox{$\halmos$}}

\newenvironment{proof}{\noindent {\bf Proof}.}{\proofbox\par\smallskip\par}

\newenvironment{remark}{\noindent {\bf Remark}.}{\par\smallskip\par}

\newenvironment{proofof}[1]{\noindent {\bf
Proof of #1}.}{\proofbox\par\smallskip\par}

\newcommand{\p}{{\mathbb{P}}}

\newcommand{\np}{{\cal G}_{n,p}}

\newcommand{\ex}{{\mathbb{E}}}

\newcommand{\al}{\alpha}

\newcommand{\ga}{\gamma}

\newcommand{\eps}{\varepsilon}
\newcommand{\la}{\lambda}

\newcommand{\co}{\mathrm{C}}
\newcommand{\out}{e^{\mathrm{out}}}
\newcommand{\td}{\mathrm{d}}
\newcommand{\D}{\bf {D}}
\newcommand{\T}{\mathrm{T}}

\newtheorem{theorem}{Theorem}[section]
\newtheorem{lemma}[theorem]{Lemma}

\newtheorem{proposition}[theorem]{Proposition}
\newtheorem{claim}[theorem]{Claim}
\title{The Evolution of the Mixing Rate}

\author{N. Fountoulakis \\
{\small School of Mathematics} \\ {\small University of Birmingham }\\
{\small United Kingdom} \\
{\footnotesize $\mathtt{nikolaos@maths.bham.ac.uk}$}
\and  B.A. Reed \\
{\small Canada Research Chair in Graph Theory,}
\\ {\small School of Computer Science,} \\
{\small McGill University, Canada}\\
{\small and} \\
{\small Laboratoire I3S, CNRS,} \\
{\small Project Mascotte, INRIA,} \\
{\small Sophia-Antipolis,} \\
{\small France } \\
{\footnotesize $\mathtt{breed@cs.mcgill.ca}$}
}

\begin{document}

\maketitle
\begin{abstract}In this paper we present a study of the mixing time of
a random walk on the largest component of a supercritical random graph, also
known as the giant component. We identify local obstructions that slow down
the random walk, when the average degree $d$ is at most 
$\sqrt{\ln~n \ln~\ln~n}$, proving that the mixing time in this case 
is $O((\ln~n/d)^2)$ asymptotically almost surely. 
As the average degree grows these become negligible
and it is the diameter of the largest component that takes over, yielding
mixing time $O(\ln~n/\ln~d)$. 
We proved these results during the 2003-04 academic year. 
Similar results but for constant $d$ were later 
proved independently by I.~Benjamini, G.~Kozma and N.~Wormald in~\cite{BKW}. 
\end{abstract}
\section{Introduction}
Given a graph $G$ with vertex set $V_n=\{1,2,..,n\}$, 
the simple random walk on $G$ is the Markov chain where 
for every edge $ij$ of $G$, the transition probability 
$p_{i,j}$ from $i$ to $j$ is $1 \over d(i)$.
I.e. we exit a vertex via a uniformly chosen edge. 
Formally, we have defined the entries of an $n$ by $n$ 
tranistion matrix  $P$ for the chain
(where $p_{i,j}=0$ if $ij$ is not an edge) and the distribution of the last 
point in a $t$ step walk from initial distribution $x_0$ is 
$x_0P^t$. 

This chain  (is ergodic and therefore) 
has a limit distribution precisely if $G$ is connected and non-bipartite.
In this case, the limit distribution $\pi$ satisfies
$\pi_i={d(i) \over 2|E(G)|}$ where $d(i)$, the degree of $i$, is the number of
edges of $G$ incident to $i$ 
(we extend this notation to sets letting $d(S)$ be the sum of the 
degrees of the vertices in $S$). Furthermore, the chain is reversible, 
as $\pi_ip_{i,j}=\pi_jp_{j,i}={1 \over 2|E(G)|}$.

We are interested in the mixing time of this chain for various random graphs on 
$V_n$. 
In this setting, we say that an event 
occurs {\em asymptotically almost surely} (a.a.s.) if its probability tends 
to 1 as $n \rightarrow \infty$. 

We consider the mixing time 
convergence with respect to the {\em total variation distance}
$d_{TV}$ between two probability distributions on $V_n$ defined as:
\[ d_{TV}\left( p^{(1)},p^{(2)} \right) = 
\max_{A\subseteq V_n} \left| p^{(1)}(A)- p^{(2)}(A)\right|. \]
Thus the mixing time of the chain is 
\[ T_{mix}(G) = \sup_{x_0} \min \{t \ : \ d_{TV}(x_0P^{t},\pi) < 1/e \}.\]
It is easy to prove that $\min \{t \ : \ d_{TV}(x_0P^t,\pi) < (2/e)^l \}
\leq l T_{mix}$.
So, $T_{mix}$ measures not only how long it takes to get to within $1/e$ of 
$\pi$, but also bounds how long it takes to get arbitrarily close to $\pi$.
Thus, it measures the rate at which the Markov chain mixes.

If we choose a graph $G$ uniformly at random from all graphs on $V_n$ then 
a.a.s. 
$$\forall ~i \in V_n~~\left|d(i)-{n \over 2}\right| \le \sqrt{n}\ln ~n$$
and
$$\forall ~i \neq j~~\left||N(i) \cap N(j)|- {n \over 4}\right| 
\le \sqrt{n}\ln~n,$$
where $N(i)$ is  the neighbourhood of vertex $i$ in $G$.
It follows easily that a.a.s. for every $j$ $\pi(j)={1 \over n}+
o({1 \over n})$ and, by counting 
the number of paths of length 2 from $i$ using the inequalities 
above, that a.a.s.  
$$\forall ~i \neq j,~~P_i^2(j)={1 \over n}+o\left({1 \over n}\right).$$
So a.a.s. $T_{mix}(G)=2$,
which is also the diameter of $G$.

We shall consider the random graph $G_{n,p}$ on $V_n$ where each edge 
is present independently with probability $p$, and hence the 
expected degree of a vertex is $p(n-1)$. In what follows we 
always couple the use of 
$p$ and $d$ by insisting that $d=pn$; this is essentially the expected degree
and makes our formulas a little easier than if we used $d$ to 
represent the actual expected degree.

In the same vein, we can prove:

\begin{theorem} \label{bignp-mix-time}
For every $p=p(n)$ with $d- \ln~n = \omega (1)$ 
we have that a.a.s. 
\[ \left|T_{mix}(G_{n,p}) -
\frac{\ln~n}{\ln~d}\right| \le 3.\]
\end{theorem}
This result improves on earlier work of Hildebrand~\cite{Hild} who had 
determined
the mixing time up to a multiplicative factor 
for $p = \Omega({\ln^2~n \over n})$. 
Its easy proof, given in a companion paper \cite{Evol2}, relies on proving inductively that 
the number of nodes at the $j$-th level of the breadth-first search tree
from an arbitrary vertex $i$ is concentrated around $d(i) d^{j-1}$ 
provided this is 
$o(n)$ (and hence generalizing the easy argument above).
\cite{Evol2} also uses similar arguments to obtain results
on the diameter of $G_{n,p}$ for $d=\omega(1)$, strengthening 
results in \cite{ChungLiu} and \cite{FernRam}.

The situation for $p\le {\ln~n \over n}$ is more problematic,
as the local structure of the graph begins to play a role. For one thing,
if $\ln~n - d = \omega(1)$ then $G_{n,p}$ almost surely has 
vertices of degree zero and hence is not connected (see for example
Theorem 7.3 in~\cite{Bol}). However, 
if $p>{1 + \epsilon \over n}$
for some $\epsilon>0$ then a.a.s. the largest component $H_{n,p}$ has 
order $\Omega(n)$ vertices 
whereas the second largest component has $O(\ln~n )$ vertices. 
So, for small $p$, we consider the simple random walk on the 
{\it giant} component $H_{n,p}$ of $G_{n,p}$. 

A second type of local structure which comes into play at this 
point are the vertices  of $H_{n,p}$ which are far away from 
any vertex of degree exceeding two. 
An easy second moment  argument, given at the end of the paper, shows 
that: 

\begin{lemma} \label{np-mix-time}
For $p<{\ln~n \over 5n}$ 
a.a.s. $H_{n,p}$ contains paths of length more than $\ln~n \over 4d$
all of whose interior vertices have degree two. 
\end{lemma}

Given such a path  $Q$
with $2l+1$ vertices, we can label these vertices with consecutive integers 
so that the midpoint $x=x_Q$ is labelled 0. Now, 
if we start our random walk at  $x$,
then we can mimic the initial part of the walk, until we touch 
an endpoint of $Q$, by the standard random walk on the integers
(where we go from $i$ to $i-1$ or $i+1$ with equal probability) starting at 0. 
It is well known (see for example~\cite{Fel} page 349) that we expect to take $l^2$ steps 
before we have 
seen an integer with absolute value $l$ in this walk. Thus if $x$ also denotes the distribution on $V_n$
that gives unit mass to $x$, it follows easily that 
$xP^t(Q)>{1 \over 2}$ for $t< \frac{1}{10}\left({\ln~n \over 8d}\right)^2$ 
and so $T_{mix}>\frac{1}{10}\left({\ln~n \over 8d}\right)^2$.

For $p\leq {\sqrt{\ln~n \ln\ln~n} \over n}$ 
this bound is larger than the diameter of $H_{n,p}$ and 
is tight up to a multiplicative factor. 
For $p$ which exceed ${\sqrt{\ln~n \ln\ln~n} \over n}$ but are
$O({\ln~n\over n})$, 
it is the diameter which correctly 
approximates $T_{mix}$ as these paths are too small to exert much influence.
So, we can determine the first order term in the mixing time precisely.

In doing so, we actually consider  a slightly different definition of
the mixing time: for $t > 0$ we let $T_t$ be uniformly
distributed on $\{0,\ldots, t-1 \}$ and we set  
$$T_{mix}'(G) = \sup_{x_0}\min \{ t: \;  d_{TV}(x_0P^{T_t},\pi) <1/e\}.$$

This is within a constant factor of  many other mixing times
(its definition is inspired by the properties of the uniform averaging
rule - see e.g.~\cite{LoW} Theorem 4.22 or~\cite{LoW1} Theorem 5.4 for further details). It is at most a constant factor larger than $T_{mix}$. 
For the modified chain in which we stay in the current state 
with probability $1/2$ and take a step with probability $1/2$ in each iteration,
$T_{mix}$ is at most a constant factor larger than this mixing time.
We have  shown that actually for the chains we are considering, a.a.s. 
$T_{mix}'$ is within a constant factor of 
$T_{mix}$ even without this modification. 
This result is a consequence of a more general result which
will appear in a separate paper~\cite{Sample}. In particular, it
answers Problem 17 in Section 4.3.3 in~\cite{AldF}.

We show:
\begin{theorem} \label{midnp-mix-time}
For every $p=p(n)$ with ${\sqrt{\ln~n \ln \ln~n}  \over n} \le 
p \le {2 \ln~n \over n}$ we have that a.a.s. 
\[ \left|T_{mix}'(G_{n,p}) -
\frac{\ln~n}{\ln~d}\right| =O\left( \left({\ln ~n \over d}\right)^2\right).\]
\end{theorem}

\begin{theorem} \label{smallnp-mix-time}
For every $\epsilon>0$,
and  $p=p(n)$ with $ {1 + \epsilon \over n}<p<{\sqrt{\ln~n \ln~\ln~n }\over n}$ 
we have that a.a.s. 
\[ T_{mix}'(G_{n,p}) =O\left(\left({\ln n \over d}\right)^2\right).\]
\end{theorem}
The proof of the first of these theorems is similar to but
more complicated than that of Theorem \ref{bignp-mix-time} and
it will be given in a companion paper \cite{Evol2}. 
In this paper,
we handle the  more delicate situation of small $p$, proving
the second theorem. 
A similar result concerning the $G_{n,m}$ model with constant average
degree was proved independently by I.~Benjamini, G.~Kozma and N.~Wormald
in~\cite{BKW}. We remark that as  lower bounds on $T_{mix}$ are also lower bounds on $T_{mix}'$ (up to a constant factor); our results are tight.
 
We close this introductory section presenting some tools which we use later in our proofs.

For any set $S$ of states of a Markov chain we define $Q(S)$ 
to be the probability we leave $S$ when we are in the steady state;
so $Q(S)= \sum_{i \in S,j \not \in S} \pi(i)p_{i,j}$ and $Q(S) \over \pi(S)$
is the probability that we leave $S$ given we are in it. 
Thus, $\pi(S) \over Q(S)$ is the expected length of a sojourn in $S$ when we
are in the steady state. It follows easily , see e.g. \cite{Cond}, 
that the mixing time of any Markov chain is at least 
$\max_{ \{S: \ 0<\pi(S) \le {1 \over 2}\}} {\pi(S) \over 10 Q(S)}.$

As in~\cite{JS}, we define the conductance of $S$, denoted $\Phi(S)$  to be 
$Q(S) \over \pi(S) \pi(V_n\setminus S)$ and the conductance 
$\Phi = \Phi(G)$ of $G$ to be $\max_{\{S: \ 0<\pi(S)<1\}} 
\Phi(S)$.
We note that for reversible chains, $\Phi(S)=\Phi(V_n\setminus S)$, so 
$$T_{mix} \ge \frac{1}{10\Phi}.$$
Moreover, M. Jerrum and 
A. Sinclair proved in~\cite{JS} that the mixing time
of an irreducible, aperiodic and reversible 
Markov chain satisfies: 
\begin{equation} \label{JerSin} 
T_{mix}' \leq \frac{C}{\Phi^2} \log \pi_{min}^{-1}, \end{equation} 
for some constant $C$. Thus, the mixing time of a reversible Markov chain is 
approximately determined by its conductance. 

In \cite{Cond}, the authors of this paper,
treading a path blazed by L. Lov\' asz and R. Kannan in~\cite{KL}, 
proved a strengthening of this result which can be used to tie down the 
mixing time of many Markov chains more precisely. 

For $p> \pi_{min}$ 
we let $\Phi(p)$ be the minimum conductance of a connected set $S$ with 
${p \over 2} \le \pi(S) \le p$ (if there is no such a set we define 
$\Phi(p)=1$).
In \cite{Cond}, we prove: 
\begin{theorem} \label{mix-time}
For an irreducible, reversible and aperiodic Markov chain on 
$V_n$ we have 
\[T_{mix}' \leq C  \sum_{j=1}^{\lceil \log \pi_{min}^{-1}\rceil} 
\Phi^{-2}(2^{-j}),\]
for some constant $C$ that does not depend on the chain. 
\end{theorem}

{\bf Remark:} The above sum can be approximated within a constant factor
by the integral $\int_{\pi_{min}/2}^{1/2}\frac{dx}{x\Phi^2(x)}$. So if
we bound $\Phi (p)$ by $\Phi$, we obtain~(\ref{JerSin}). However, the
bound of Theorem~\ref{mix-time} is often tighter. Such is the case for
the chain considered in this paper. 
\vskip0.2cm

We apply Theorem \ref{mix-time} to deal with the mixing time for 
small $p$.
To do so, we need 
to bound  from above the conductance of connected sets of various sizes. 

This completes our preliminary remarks. In the next section, we 
discuss the precise results on the conductance of $H_{n,p}$ that 
we need to prove Theorem \ref{smallnp-mix-time} and show that they do indeed imply 
this result. In Sections~\ref{SimpleFacts}-\ref{ProofCore}, we prove these results.  
In Section~\ref{LongPaths}, we prove Lemma~\ref{np-mix-time} which shows that the bounds
of Theorem \ref{smallnp-mix-time} are tight.
We close this section with the statement of Talagrand's inequality 
(see for example inequality (2.43) p. 42 in~\cite{Rg}), 
which we will use at various points in our proofs to derive
concentration bounds: there exists a constant $\gamma > 0$ such that 
for any $t>0$, if $X$ is a binomially distributed random variable:  
\begin{equation} \label{Talagrand}
\p [|X -\ex [X]|>t] \leq 4 e^{-\gamma t^2/(\ex[X]+t)}.
\end{equation}

\section{The Evolution of the Conductance}

We focus now on the conductance of the connected subsets of $H_{n,p}$. 
To this end, we let $e^*$ be the number of edges of $H_{n,p}$ and recall that 
for non-biparite $H_{n,p}$,
for every vertex  $v$, $\pi(v)={d(v) \over 2e^*}$. 

Letting 
$e(S) = |\{\{i,j\} \in E:\ i,j \in S\}|$ and
$\out(S)=|\{\{i,j\} :\ i \in S, j \not\in S\}|$,
we have that if $H_{n,p}$ is non-bipartite: $\pi(S)=
{d(S) \over 2e^*}= {\out(S)+2e(S) \over 2e^*}$ and
$Q(S)={\out(S)  \over 2e^*}$.
Thus for such non-bipartite $H_{n,p}$, $\Phi(S)
= {\out(S) \over (2e(S)+\out(S))\pi(V_n\setminus S)}$
which is within a factor of 2 of ${\out(S) \over 2e(S)+\out(S)}$ 
and so to bound conductance  we need to know about the behaviour of 
these two variables. 

The advantages of focussing on these variables individually, rather than 
on conductance itself are two-fold. The first is that they are defined for 
all $H_{n,p}$ not just for non-bipartite ones. The second is that we can 
quickly see that it is $e^{out}(S)$ which gives us difficulty. Indeed, 
standard concentration results easily yield: 

\begin{lemma} \label{B*}
There exists an absolute constant $l$ such that for any $p$
with ${1 \over n} <p<{2\ln~n \over n}$, 
a.a.s. every connected set $S$ in $G_{n,p}$ 
satisfies $e(S)\leq ld|S|$. Furthermore, 
if $|S|=o((\ln~n)^2)$ then  $e(s) \le 
2|S|$.
\end{lemma}

The easy proof of Lemma \ref{B*} is given in the next section. 
The behaviour of $e^{out}(S)$ is more problematic. However, more 
careful counting arguments allow us to show:

\begin{lemma} \label{E} 
There exists constants
$\epsilon>0, c>0$ and $d_0> 1$ such that for 
every $p=p(n)$ satisfying
${d_0 \over n} \leq p \leq {2\ln~n \over n}$, 
a.a.s. every connected subset $S$ 
of $H_{n,p}$ with $|S|\geq \frac{c\ln n}{d}$, 
and $d (S)\leq \frac{d(H_{n,p})}{2}$ satisfies $\out(S)\geq \epsilon d|S|$. 
\end{lemma}

\begin{lemma} \label{F}
For any two constants $c_1$ and $c_2$ 
with $1<c_1<c_2$ 
there exist $\epsilon, A>0$ such that
a.a.s. for $d=d(n)$ 
lying between $c_1$ and $c_2$,
every connected subset $S$ of $H_{n,p}$ with 
$|S|\geq A{\ln n}$ and $d(S)\leq \frac{d(H_{n,p})}{2}$ satisfies
$\out(S)\geq \epsilon |S|$. 
\end{lemma}

Since $H_{n,p}$ is connected  we know that $e^{out}(S)$ is at least one for 
all strict subsets of $V(H_{n,p})$.
Armed with this fact and  Lemmas \ref{B*},\ref{E}, \ref{F}, 
we can now give the:
\vskip0.2cm

\begin{proofof}{\bf  Theorem \ref{smallnp-mix-time}} 
It is well known that $H_{n,p}$ is a.a.s
non-bipartite.
So we assume this to be the case
and apply Theorem  \ref{mix-time}, to obtain:

\[T_{mix}'(H_{n,p}) \leq C  \sum_{j=1}^{\lceil \log \pi_{min}^{-1}\rceil} 
\Phi^{-2}(2^{-j}).\]
We now bound the sum $\sum_{j=1}^{\lceil \log \pi_{min}^{-1}\rceil} 
\Phi^{-2}(2^{-j})$.

Applying Lemma \ref{B*} with $S=V$, we have that $e^*=O(dn)$.
Since $H_{n,p}$ is connected, $\pi_{min}$ is at least
${1 \over e^*}=\Omega({1 \over dn})$.
So $\log~\pi_{min}^{-1}=O(\ln~n)$.
We claim that we can apply 
Lemmas \ref{B*},\ref{E},
and \ref{F} to show that there exists 
absolute constants  
$r$ and $c$  such that if $\Phi^{-1}(S)>r$ then 
$\pi(S) \le {c \ln~n \over d^2n}$.
So letting $I$ be the set of $j$
such that $2^{-j}$ lies between 
$\pi_{min}$ and $c \ln~n \over d^2n$,
our sum is bounded by $O(\ln~n)+
\sum_{j \in I} \Phi^{-2}(2^{-j})$.

But, $\Phi(S) \ge {1 \over e^*\pi(S)}$, for such $S$ because $H_{n,p}$ is connected.
So, since $e^*=O(dn)$,
$\Phi^{-2}(2^{-j})=O(d^2n^22^{-2j})$ uniformly
for $j \in I$. Since this
is a geometric sum, it is of the 
same order as its largest term,
which is $O\left(\left({\ln~n \over d}\right)^2\right)$,
as required. 

It remains to prove our claim.

Recalling how we have expressed $\Phi$ in terms of $e^{out}$ and $e$, and 
combining Lemmas \ref{B*}, \ref{E}, and \ref{F}
for $c_1={1+\epsilon \over n},c_2={d_0}$ 
we have:

For any $\epsilon>0$,
there is a $c$ and an $r$
such that for any $p=p(n)$ which is 
at most ${\ln~n \over n}$ and at least 
$ {1 + \epsilon \over n }$, a.a.s
any connected set $S$ with at least  
$\frac{c \ln n}{d}$  vertices  and
$\pi(S) \le {1 \over 2}$ satisfies
$\Phi^{-1}(S)\le r$ (note that when $d$
is less than $d_0$, it can be treated as a constant).
 
But 
since $\Phi(S)>{1 \over 3}$ if
$e^{out}(S)> e(S)$, the second half of Lemma \ref{B*} tells us that a.a.s. every
such small set whose conductance is small
satisfies 
$d(S)\le 4|S|$ and so $\pi(S) =
O\left({|S| \over dn}\right)$.
Combining this with our bound on the size of
S proves the claim and the lemma.
\end{proofof}

\vskip0.5cm
Lemma \ref{E} is easy to prove using Talagrand's inequality.
Lemma \ref{F} is much more difficult. To prove it we consider the {\em core}
of $G_{n,p}$ which is the maximal subgraph of $G_{n,p}$ all of whose vertices 
have degree at least 2. We prove some results about the expansion properties
of the core and then translate them into results on $H_{n,p}$ which imply 
Lemma \ref{F}. In doing so, we condition on the degree sequence of the 
core and we use the Bender-Canfield model 
for graphs with a given degree sequence. We discuss this model in
Section~\ref{core-analysis} but readers may consult~\cite{BenC} for
further details.

In the next section we present some well-known properties of $G_{n,p}$ and
prove Lemmas \ref{B*}, Lemma  \ref{E} and the following, of which we shall 
have need.  
\begin{lemma} \label{total-degree}
For every fixed $d>1$ there exists $L >0$ such that if $p= {d \over n}$ then a.a.s. 
every connected set $S \subseteq V_n$ in $G_{n,p}$ with $n \geq |S|\geq  \ln n$ 
has $\td (S) \leq L|S|$.
\end{lemma}
The proof of Lemma \ref{E} is also presented in Section \ref{SimpleFacts}.

Now Lemma \ref{np-mix-time} guarantees that a.a.s. $G_{n,p}$ contains
paths which are connected sets with conductance 
between $\frac{1}{\lfloor\frac{\ln n}{d} \rfloor}$ and 
$\frac{2}{\lfloor\frac{\ln n}{d} \rfloor}$.  
Since $\frac{x}{m+x}$ us non-decreasing in $x$, Lemma \ref{A} below tells
us that every set $S$ of order $k\leq A_p\frac{\ln n}{d}$ has conductance at
least $\frac{1}{2k+1}=\Omega \left(\frac{d}{\ln~n} \right)$.
On the other hand Lemmas \ref{E} and \ref{F} tell us that the minimum 
conductance a set of order at least $A_p \ln n /d$ can have is a.a.s. $\Omega (1)$. 
So, if $\Phi (H_{n,p})$ denotes the minimum conductance of a subset of
$H_{n,p}$ we have:
\begin{theorem} Whenever $d=d(n)>1+\Theta(1)$, a.a.s,  
\[\Phi (H_{n,p})=\Theta \left( \min\left\{\frac{d}{\ln n},1\right\} 
\right).\]
\end{theorem}

\section{Some Simple Facts about $G_{n,p}$} \label{SimpleFacts}

To prove Lemma \ref{B*}, we simply combine the following two lemmas. 

\begin{lemma} \label{A} 
For every $p=p(n)$ with ${1\over n}< p \le {2\ln~n  \over n}$,  
a.a.s. every connected set $S\subseteq V_n$  of $G_{n,p}$ with 
$|S|\leq \frac{n}{60d^2}$ satisfies $|S|-1 \leq e(S)\leq 2|S|$. 
\end{lemma}

\begin{proof} 
The expected number of sets of $2k$ edges spanning a set of
$k$ vertices is 
${n \choose k}{{k \choose 2} \choose 2k}p^{2k}\leq 
\left(\frac{ne}{k}\right)^k \frac{e^{2k}k^{4k}}{k^{2k}} \frac{d^{2k}}{n^{2k}}=
\left(\frac{e^3kd^2}{n}\right)^{k}$. For $k \le \frac{n}{60d^2}$ this is 
less than $2^{-k}$. For $k$ less than  $\sqrt{n} \over e^3d^2$, it is 
less than $\sqrt{n}^{-k}$. 
The result follows by the first moment method.
\end{proof}

\begin{lemma} \label{B}
There exists a constant $l$ such that 
for any $d=d(n)$ between $1$ and $2\ln~n$, 
and $p={d \over n}$, 
a.a.s. every set $S$ with $|S|\geq \frac{ n}{60d^2}$ satisfies $e(S)\leq ld|S|$.
\end{lemma}
\begin{proof} 
If $e_k$ denotes the number of edges on a set of $k$ vertices,
then $\ex[e_k]=p{k\choose 2}\leq\frac{dk}{2}$. 
As $e_k$ is binomially distributed,
we can use Talagrand's inequality (\ref{Talagrand}) 
to show that 
\[\p [e_k >\ex [e_k] + t] \leq 4\exp 
\left(-\frac{\ga t^2}{\ex [e_k]+t} \right), \]
for some universal constant $\ga >0$.   
In particular, for a natural number $l$ which will be determined soon  
\begin{eqnarray*}
\p [e_k > l dk]\leq \p [e_k >\ex [e_k] + (l-1/2)dk] 
\leq 4 \exp \left( -\frac{\ga (l-1/2)^2 (dk)^2}{ldk}\right) 
\leq 4 \exp \left(-\frac{\gamma ldk}{2} \right),
\end{eqnarray*}
for $l$ sufficiently large. 
Therefore the expected number of sets with $k$ vertices and at least 
$ldk$ edges is  at most $4{n \choose k}e^{-\frac{\ga ldk}{2}}\leq 
4\left(\frac{ne}{k} e^{-\frac{\ga l d}{2}} \right)^k \leq 
4\left(60 ed^2 e^{-\frac{\ga ld}{2}} \right)^k$ and 
\[4\sum_{k\geq \frac{\ln n}{cd^2}} \left(60 ed^2 e^{-\frac{\ga ld}{2}}\right)^k
=o(1), \]
for $l$  a sufficiently large constant.
\end{proof}

We will need the following simple fact in the proof of Lemma \ref{E}.
It will also be useful later. It is an immediate consequence 
of, for example, Talagrand's inequality.

\begin{lemma}
\label{honeybunny}
For $p>{1 \over n}$, a.a.s $G_{n,p}$ has 
$({1 \over 2}+o(1))pn^2$ edges.
\end{lemma}

\begin{proofof} {\bf Lemma {\ref{E}}} 
%
We specify our choice of $\epsilon < {1 \over 6}$ below.
So we can assume $e^{out}(S) =e^{out}(V-S)$ is at most 
$\frac{dn}{6}$ as otherwise we are done.
Since $d(S) \le {d(H_{n,p}) \over 2}$,  by Lemma \ref{honeybunny},
$d(V-S) \ge \frac{dn}{3}$ and so $e(V-S) \ge \frac{dn}{12}$. So, we 
are done by Lemma \ref{B*} for $S$ with $|S| \ge n -\frac{n}{12l}$. 

Key to our analysis is the fact that there are $k^{k-2}$ labelled trees on 
$k$ vertices. Thus the expected number of trees of $\np$ with $k$ vertices is 
\begin{equation} \label{tree-exp}
{n \choose k} k^{k-2} \left(\frac{d}{n} \right)^{k-1} \leq n (ed)^k.
\end{equation}
Having fixed both a set $S$ of $k \le n -\frac{n}{12l}$ 
vertices and the set of edges
within $S$, the expected value of $\out (S)$ is $dk(1-k/n)$ which is at least 
$dk/12l$. Since $\out (S)$ is binomially distributed,
Talagrand's inequality (\ref{Talagrand}) yields 
\[\p [\out(S)<\ex [\out(S)] - t] \leq 4\exp 
\left(-\frac{\ga t^2}{\ex [\out(S)]+t} \right). \] 
Setting $t=dk/24l$, we obtain: 
\begin{equation} \label{Tal}
\p \left[\out (S)< \frac{dk}{24l}\right] \leq 4\exp \left(-\ga \frac{dk}{100l} 
\right). 
\end{equation}
Combining (\ref{tree-exp}) and (\ref{Tal}), we obtain that 
for $k \leq n(1-1/12l)$, 
the expected number of connected sets $S$ of $G_{n,d/n}$ with $|S|=k$ and 
$\out (S)\leq dk/24l$ is at most 
\[4n (ed)^k e^{-\ga dk/100l}. \]
If we choose $d_0$ such that $ed e^{-\ga d/100l}< e^{-\gamma d/200l}$
for every $d>d_0$, and set
$c=\frac{600l}{\ga}$ then 
\[ 4n \sum_{k \geq \frac{c\ln n}{d}} (e^{1-\ga d/100l}d)^k \le
4n \sum_{k \geq \frac{c\ln n}{d}} (e^{-\gamma d/200l})^k =
4n \sum_{k \geq \frac{c\ln n}{d}} O\left( \frac{1}{n^3}\right)=
o(1).\] 
\end{proofof}
\begin{proofof} {\bf Lemma {\ref{total-degree}}}
Since $\td (S)=2e(S)+\out(S)$ and we have an upper bound on $e(S)$ by
Lemmas \ref{A} and \ref{B}, it suffices to bound $\out (S)$. 
For a set $S$ with $k$ vertices $\ex[\out(S)]=pk(n-k)=dk\left(1-\frac{k}{n}
\right)\leq dk$. By
Talagrand's inequality (\ref{Talagrand})
\begin{eqnarray*} 
\p [\out(S)>l dk]&\leq & \p [\out(S)>\ex [\out(S)] + (l-1)dk] \leq 4 \exp
\left(- \frac{\ga((l-1)dk)^2}{ldk} \right) \\
&\leq &
4 \exp\left(- \frac{\ga ldk}{4} \right). 
\end{eqnarray*}
Hence the expected number of connected sets $S$ with $k$ vertices and 
$\out(S)$ which is at least $ldk$ is at most
\[4n (ede^{-\frac{\ga ld}{2}})^k. \]
Choosing $l>{20 \over \ga}$ so that $ede^{-\ga l d}< e^{-3}$ we 
obtain:
\[4n\sum_{k= \lceil \ln n\rceil}^{n} (ede^{-\frac{\ga ld}{2}})^k \leq 
4n^2 e^{-3\ln n} =o(1), \]
and this concludes the proof of the lemma.
\end{proofof}

\section{ The Core of the Proof} \label{ProofCore}

In this section we obtain a lower bound on $\out(S)$ for any sufficiently      
large connected subset of $H_{n,p}$, for any $d$ bounded between two constants.
To do so, we mainly investigate 
the expanding properties of the core of the biggest component of a 
$G_{n,p}$ random graph, which we will denote by $\co (H_{n,p})$. 
This is the maximal subgraph of $H_{n,p}$ having minimum degree at 
least 2.
We let $N=\left|V(\co (H_{n,p})) \right|$ and $M=\left|E(\co (H_{n,p})) 
\right|$ and recall (see \cite{Pit}) that 
a.a.s. both $M$ and $N$ are $\Theta(n)$, for the $d$ within the range we are 
interested. 
The following lemma bounds below $\out (S)$ for   
any connected $S$ sufficiently large which is a subset of the vertex-set 
of the core. 
\begin{lemma} \label{e-out}
For every $d=d(n)$ with 
$1<c_1<d<c_2$ for two constants $c_1$ and $c_2$,
and $\la \in (0,1)$
there exist constants 
$\al_0 = \al_0(d, \la)>0$ and $C_1 = C_1 (d,\la)$ such
that a.a.s. every connected $S \subseteq V(\co (H_{n,p}))$ with  
$C_1 \ln n\leq |S|\leq \la N$ is joined to 
$\co(H_{n,p})-S$ by at least $ \al_0 |S|$ edges.
\end{lemma}

The proof of this lemma is postponed until the next subsection. 
It is reasonably easy because the fact that the core has minimum
degree two gives it strong expansion properties.

We are interested in the expanding properties, of
connected subsets of $H_{n,p}$ rather than its core. 
That is we 
would like to bound below $\out (S)$ in the case where $S$ is an arbitrary 
sufficiently large connected subset of $V(H_{n,p})$, using Lemma \ref{e-out}.
To this end, we note that
the components of $H_{n,p}- \co (H_{n,p})$ are
trees. Each such tree has a unique vertex that is adjacent to a vertex of the 
core; we say that the tree is {\em rooted} at that specific vertex. 
We call these trees {\it decorations} of the core.
To apply  the above lemma 
we need to prove that any sufficiently large connected subset
of vertices of $H_{n,p}$ has at least a certain proportion of its vertices 
belonging to the core. 

We show that this is indeed the case 
in the following lemma,
which bounds the probability that a sufficiently big connected set of vertices 
of the core has at least $l$ times more vertices in the decorations dangling
from its vertices.

\begin{lemma}\label{halo} For every fixed $d=d(n)$ 
with $1<c_1<d<c_2$ for some constants $c_1$ and $c_2$,
there are constants $\chi, l >1$ 
such that a.a.s. every connected subset $S$ of 
$H_{n,p}$ with $\chi \ln n\leq |S|\leq N$ is such
that the number of vertices belonging to $\co (H_{n,p} )$ is at least $|S|/l$.
\end{lemma}
\begin{proof}
Clearly, if we choose a counterexample $S$ so as to minimize the number 
of decorations which it intersects but does not contain, 
then it partially contains at most one decoration. 
I.e., it suffices to prove that for  some $\chi>0$,
for any sufficiently large integer $t$, 
the expected number of trees of $G_{n,p}$ 
with $tk \geq \chi \ln n$ vertices, 
such that $(t-1)k$ of these vertices are incident to no edges off the tree
and induce a forest tends to 0 as $n\rightarrow \infty$.    
We can see that the expected number of such trees is bounded above by 
\begin{equation} \label{ExpBadTrees}
{n \choose tk} {tk \choose k}
(tk)^{tk-2}p^{tk-1}(1-p)^{(n-tk)(tk-k)+{tk-k \choose 2} - (tk-k)}.
\end{equation}
If $tk < \epsilon n$, then 
using (\ref{tree-exp})
we can see that for any given $t$ and $k$ this is bounded above by 
$$nd^{tk}e^{tk}{tk \choose k}(1-p)^{ntk(1-1/t)(1-\epsilon)}.$$
This is at most 
$nz_d^{tk}{tk \choose k}e^{c_2(1-\epsilon)/t}$ for 
$z_d=e^{1+\ln~d-d(1-\epsilon)}$. We note that if $\epsilon$ is chosen
to be sufficiently small, then $z_d < z_{c_1}<1$. 
By making  $t$ large enough 
we can make 
${tk \choose k}e^{c_2(1-\epsilon)/t}$  
smaller than $(1/z_{c_1})^{tk/2}$.
Thus, we see that the expected number of such trees having at most 
$\epsilon n$ vertices is at most 
$$n z_{c_1}^{\chi \ln n /2}.$$
Setting $\chi = 4/\ln(1/z_{c_1})$ clearly suffices.

If $tk \ge \epsilon n$, we need to use more accurate bounds. Since $S$
lies in $H_{n,p}$ its number of vertices is a.a.s. at most $\la n$ for
some $\la = \la (d)<1$.  
Then writing $tk=\alpha n$,  (\ref{ExpBadTrees}) is bounded above by 
$$O(n)\left({ 1 \over \al^{\al} (1-\al)^{1-\al}} \right)^n {tk \choose
k} 
(\al n)^{\al n}
\left({d \over n} \right)^{\al n} e^{-nd\al (1-\al)(1-1/t) - n d\al^2/2 + n3d \al^2/2t
},$$
uniformly for any $\al \in [\epsilon,\lambda]$. 
In turn, this is bounded above by 
$$O(n){\al n \choose \al n/t} e^{n(h(a,d)+3c_2/2t)},$$
where $h(a,d)=-(1-\al)\ln (1-\al)+\al \ln~d -d\al(1-\al)-d\al^2/2$.
Elementary calculations show that for each $d>1$ there exists $c_d>0$ such that 
$h(a,d)< - c_d $, for every $\al \in [\epsilon,\lambda]$. 
Taking $t$ sufficiently large so that 
${\al n \choose \al n/t}e^{n3c_2/2t} < e^{nc_d/2}$,
it turns out that (\ref{ExpBadTrees}) is no more than 
$$O(n)e^{-nc_d/2},$$
whenever $\eps n \leq tk \le \la n$.
The result follows.
\end{proof}

In order to prove Lemma \ref{F}, we will need another result of a similar 
flavour.

A {\it dangling} tree is a rooted tree $T$ of $G_{n,p}$ all of whose non-root nodes
are incident only to the edges of $T$.

\begin{lemma}\label{hola} For every fixed $d=d(n)$ 
with $1<c_1<d<c_2$ for some constants $c_1$ and $c_2$ there 
is a constant $I$ such that setting $p={d \over n}$,
a.a.s. for every $i \ge I$, the number of nodes of 
$G_{n,p}$ in dangling trees of size $i$ or greater is
less than $n \over i^3$. 
\end{lemma}
\begin{proof}
Using (\ref{tree-exp})
we can see that for any given $j=O(\ln~n)$, the expected number of 
dangling trees of size $j$ is  bounded above by:
$$nd^{j}e^{j}(1-p)^{(j-1)(n-j)}=O(nd^je^je^{-dj}).$$
As in the proof of the last lemma, this is falling 
exponentially quickly as $j$ increases. 
It follows that a.a.s every dangling tree 
has size $O(\ln~n)$.
Furthermore, a straightforward second moment calculation shows that
assymptotically almost surely for every $j$, the number
of dangling trees of size $j$ is at most twice its expected value. Provided 
this condition holds, so does the conclusion of the lemma.
\end{proof}

With these three results in hand, we can give the 
\vskip0.2cm

\begin{proofof}{\bf Lemma \ref{F}}
To obtain a lower bound on $\out (S)$ 
where $S$ is a sufficiently 
large subset of $V(H_{n,p})$, we want
to use Lemma \ref{halo} to argue that at least 
$|S|/l$ of its vertices belong to the core and then apply  Lemma \ref{e-out}
to deduce that $S$ has at least $\al_0 |S|/l$ edges to $\overline{S}$. 
Of course such an approach is valid whenever $S$ contains a proportion of 
vertices of the core that is bounded away from 1. As we are interested in 
bounding the conductance of sets having $\pi (S)\leq 1/2$, it might be the
case that $S$ contains almost all of the core or even all of it. In particular,
this might be the case when $d$ is close to 1 (see \cite{Pit} for precise 
bounds on the proportion of vertices of $H_{n,p}$ belonging to the core). 
Hence, in this case we must argue in a different way as we shall see below. 

First we need the following:
\begin{claim} \label{TreeExp}
Suppose that $p$ is as in Lemma~\ref{halo}. 
Let $T$ denote a collection of vertex-disjoint trees which are
induced in $H_{n,p}$ and let $t$ be the total number of vertices they
contain. There exist constants $\chi, l>0$, such that every such $T$ with all
its vertices, except one per tree, incident only to edges
within $T$ and $t \geq \chi \ln~n$ has a.a.s. $|T|>t/l$. 
\end{claim}
\begin{proof}
The proof of this claim is essentially the same as that of
Lemma~\ref{halo}, except that in (\ref{ExpBadTrees}) we replace 
the factor $(tk)^{tk-2}$ by $k(tk)^{tk-k-1}$, which is the number of
forests we can build on $tk$ vertices, with $k$ particular vertices
belonging to different trees. We omit the details. 
\end{proof}

By the above remarks,  to complete the proof of Lemma \ref{F} 
we need only prove for $S$ satisfying  
$|S \cap V(\co (H_{n,p}))| > (1-\tau_0)N$, 
for a constant $\tau_0=\tau_0(d)$ to be specified later.

We note that since $S$ is connected,
for any vertex $w$ in $V(\co (H_{n,p})-S)$ the union of the decorations 
attached at $w$ is disjoint from $S$. Let $e_C (\overline{S})$ denote
the number of edges in  
$V(\co (H_{n,p})-S)$ along with the attached decorations.

Furthermore, any component
of $H_{n,p}-S$ disjoint from the core is a dangling tree rooted at a 
vertex which has a neighbour in $S$. Let $t_S$ be the total number of
vertices involved in these trees. If $T_S$ is the number of these
trees, then $e^{out}(S)\geq T_S$. We show that $T_S \ge
e^*/4l$, where $l$ is the constant that appears in the above claim.   

Indeed, since $\pi (V(H_{n,p})-S)\geq 1/2$, we have
\begin{equation} \label{VertTrees}
e^*\leq d(V(H_{n,p})-S) \leq 2t_S + e^{out}(S)+2e_C(\overline{S}).
\end{equation}
We may assume that $e^{out}(S)\leq e^*/4$, as otherwise we are done. 
Since $V(\co(H_{n,p})-S)$ contains at most $\tau_0 N$ vertices, 
Lemma~\ref{hola} implies that $2e_C(\overline{S})\leq e^*/4$, provided
that we choose $\tau_0$ small enough. Plugging these bounds into 
(\ref{VertTrees}), we deduce that $t_S \ge e^*/4$.  
Claim~\ref{TreeExp} yields the bound on $T_S$ and the proof of the
lemma is now complete.  
\end{proofof}

\subsection{The expanding properties of $\co (H_{n,p})$} 
\label{core-analysis}

In this section we investigate the expanding properties of the core of 
the giant component of a $G_{n,d/n}$ random graph, where $d$ is bounded 
between two constants. 
We are aiming towards the proof of Lemma \ref{e-out}.

We actually consider the core of $G_{n,p}$ rather than
the core of $H_{n,p}$  because 
it is known that $\co (\np)$ conditioned on its degree sequence is uniformly
distributed over all graphs having this degree sequence. 
This follows from
Proposition 2.1(b) in \cite{Psw} that describes the distribution of the graph 
that remains from the stripping-off process conditioned on its degree sequence. It is also known (again see
\cite{Psw}) that a.a.s.
the difference between these two graphs consists 
of a set of $O(\ln~n)$ cycles each with at most $O(\ln~n)$ 
vertices. So it is enough to prove the variant of Lemma~\ref{e-out} obtained by replacing the core of 
$H_{n,p}$ by the core of $G_{n,p}$. 

In doing so, we will use the configuration model 
of Bender and Canfield (see
\cite{BenC}).  Suppose we want to 
analyze a uniformly random graph 
$G_N$ on $V_N = \{1,2,...,N\}$ with a given degree sequence
$\{d_1,...,d_N\}$. 
For each $i \in V_N$ we take  $d_i$ copies of $i$, thus forming a set $P_N$  of $2M$ points. 
A perfect matching on $P_N$ corresponds to a multigraph on $V_N$ where  an edge between copies of $i$ and $j$
yields an edge from $i$ to $j$. Note that this may 
create loops or multiple edges.  
We consider the random multigraph $G_N'$ that comes out of 
a uniformly random perfect matching on $P_N$.
As shown by McKay and Wormald in~\cite{McW}, under certain
technical conditions,  results proven in this
model can be transferred back to a uniformly chosen simple
graph on the given degree sequence.

The key to our lemma is the following, which 
is very easy to prove by considering this correspondence.
 
\begin{proposition} \label{no-edges-across}
For each integer $N\geq 1$ let $(d_1,\ldots , d_N)$ 
be a degree sequence on the set $V_N=\{1,\ldots, N \}$ and let 
$G_N$ be a random graph having this degree sequence.
Assume that 
\begin{enumerate}
\item For each $N$ and for $i=1,\dots , N$ we have $d_i \geq 2$. 
\item 
$2M=\sum_{i=1}^{N}d_i \leq \sum_{i=1}^{N}d_i^2 \leq C N$, for some $C>0$. 
\item $\max \{ d_i\}_{1\leq i \leq N} \leq \ln N$, for $N$ sufficiently large. 
\end{enumerate}
Then uniformly for every set of vertices $S$ with 
$\td (S)=\sum_{i \in S}d_i$ being an even 
number, we have 
\[\p [\mbox{There are no edges between $S$,$V_N\setminus S$ in $G_N$}] 
=O\left( {M \choose \td (S)/2}^{-1} \right). \]
(In fact the constant involved in the $O(1)$ term depends only on $C$.)
\end{proposition}
\begin{proof} 
Given a degree sequence and a set $S$
of vertices  satisfying the conditions of the lemma,
we generate the perfect matching on  $P_N$, which yields
$G_N'$, one edge at a time.
 
Exposing the edge out of the first vertex of $S$, we see that the probability
it is joined to a vertex of $S$ is $\frac{\td (S)-1}{\td (V(\co(\np)))-1}$. 
Given 
that this edge stays within $S$, we take a third vertex of $S$ and note that 
the probability that it is joined to a vertex in $S$ is 
$\frac{\td (S)-3}{\td ( V(\co (\np)))-3}$. More generally:

\begin{eqnarray*} 
\lefteqn{
\p [\mbox{There are no edges between $S$, 
$V_N\setminus S$ in $G_N'$}]=} \\
& &\frac{\td(S)-1}{2M-1} \ \frac{\td(S)-3}{2M-3} \cdots
\frac{1}{2M-\td(S)+1}  \\
&=& \frac{\td(S)!}{(\td (S)/2)!2^{\td(S)/2}} \ 
\frac{(2M-\td(S))!}{(M-\td(S)/2)!2^{M-\td(S)/2}} \
\frac{M!2^M}{(2M)!} \nonumber \\
&=& {M \choose \td(S)/2} / {2M \choose \td(S)} \leq {M \choose \td(S)/2}^{-1}. 
\end{eqnarray*} 

Assumptions $(2)$ and $(3)$ along with the main theorem in \cite{McW} 
transfer this bound to the space of a $G_N$ random graph,
and we are done.
\end{proof}

For each $n\geq 1$, and any $C,c,  \varepsilon >0, B$, 
we define $E_n=E_n(C,c, \varepsilon)$ to be the set of
graphs on $V_n$ such that their core is non-empty, it has $N$
vertices and $M$ edges, where
$N \ge cn$ and 
$N/M \leq {1 - \epsilon}$, maximum
degree at most $\ln N$ and moreover, the sum of the squares of the 
degrees in the core is no more than $CN$ (see condition (2) in the
previous proposition). To prove Lemma \ref{e-out} we will condition on
the event $E_n$ for a specific choice of $C, c$ and $\varepsilon$: \\ 
\begin{remark}
There exist $C, c, \varepsilon$ such that for every $1<c_1 < d < c_2$  
the event $E_n$ occurs a.a.s. (see~\cite{Bol},~\cite{Pit}).
\end{remark}

Now we are ready to proceed with the proof of Lemma \ref{e-out}:

\begin{proofof}{\bf Lemma \ref{e-out}}
In fact 
we shall 
prove that every $S \subseteq V(\co (\np))$ 
with $G[S]$ being connected and $\la N\geq |S|\geq C_1 \ln n$ 
has more than $\al_0 |S|$ edges joining it with 
$V(\co (\np))\setminus S$, for some $\al_0$, $C_1$ which will be specified  
during our proof. For the time being we assume that 
$C_1 \geq 1$,
with the prospect of using Lemma \ref{total-degree}.  

Now for any integer $s$ between $C_1\ln n$ and $\la N$, let 
$X_s = X_s (d, \al_0)$
be the number of subsets $S$ of $V(\co (\np))$ with $|S|=s$, $d(S)\leq Ls$, 
$e(S)\geq s-1$ and having at most 
$\al_0 |S|$ edges joining it with $V(\co (\np))
\setminus S$. Considering such sets is sufficient, since 
by Lemma \ref{total-degree} a.a.s. every connected set of vertices $S$ with 
$|S|\geq C_1 \ln n$ has total degree no more than $L|S|$. 
We will condition on the event $E_n$ and, more
specifically, we will show that $\ex [X_s\ | \ E_n]=o(1/n)$ 
uniformly for every $s$ between $C_1\ln n$ and $\la N$. 
To do so, we shall estimate the conditional expectation of $X_s$ given
a degree sequence of the core of a graph in $E_n$ (in
the conditional probability space of the event $E_n$), which we denote
by $\tilde{\ex}[X_s]$. 
Proving that $\tilde{\ex} [X_s]=o(1/n)$, i.e. that the random variable 
$\tilde{\ex}[X_s]$ is bounded above by a function that is $o(1/n)$, 
uniformly for any degree sequence of the core of a graph in $E_n$ 
will be sufficient. Then the lemma will follow since $E_n$ occurs a.a.s., if
we choose $C, c$ and $\varepsilon$ as in the above remark.

We proceed with the estimation of $\tilde{\ex} [X_s]$, where 
$C_1\ln n \leq s \leq \la N$. For any  $\al \leq \al_0$,
if the set $S$ has total degree $t$, then we can choose 
the edges that will be adjacent to $V(\co (\np ))\setminus S$ in at most 
${t \choose \al s}$ ways. Having specified these elements of the total degree, 
now the available total degree in $S$ is equal to $t'=t-\al s$.
Then the probability (i.e. the conditional expectation of the
indicator random variable that is equal to 1 on the event) that every
other edge is not adjacent to $V(\co (\np))\setminus S$ is 
by  Proposition \ref{no-edges-across} 
$O\left({M \choose t'/2}^{-1}\right)$, uniformly over all sets $S$ as above.
This is the case because the premises of Proposition \ref{no-edges-across} 
are satisfied, by the definition of $E_n$.  

Since $s-1 \leq t'/2 \leq M-N+s$, the above bound is 
$O(1)\max\left\{ {M \choose s-1}^{-1}, {M \choose N-s}^{-1} \right\}$.
Note that by Stirling's formula we have 
\begin{eqnarray} \label{bound1}
\frac{{N \choose s}}{{M \choose s-1}} &=& 
\frac{(s-1)!}{s!} \ \frac{N!}{M!} \ \frac{(M-s+1)!}{(N-s)!} 
= O(1) \frac{1}{s} \frac{N^N}{M^M} \ \frac{(M-(s-1))^{M-(s-1)}}
{(N-s)^{N-s}} \nonumber \\
&=& O(1) \frac{1}{s} \ N^s \ \left(1-\frac{s}{N}\right)^{-N+s} 
\frac{M^{M-s+1}}{M^M} \left(1-\frac{s-1}{M} \right)^{M-(s-1)} \nonumber \\
&=& O(1) \frac{M}{s} \ \left( \frac{N}{M} \right)^s \exp 
\left(s -\frac{s^2}{N} - s + 1 + \frac{s^2}{M} \right) 
= O(1) \frac{M}{s} \ \left( \frac{N}{M} \right)^s.
\end{eqnarray}
Since ${N \choose s}={N \choose N-s}$ and ${M \choose s}={M-s+1 \over s} {M \choose s-1}$, applying (\ref{bound1}) with $N-s$ in place 
of $s$ yields:
$$
 \frac{{N \choose s}}{{M \choose N-s}}
=O(1)\frac{M}{M-N+s+1} \left(\frac{N}{M}\right)^{N-s}.$$

Provided $\frac{N}{M} \le (1-\epsilon)$,
both these bounds are $O(Me^{-\beta s})$ for some 
$\beta=\beta(\epsilon, \lambda)$.

We are now ready to bound $\tilde{\ex} [X_s]$, for any $s$ as above.
Let $\D = \D (s, \al_0)
=\{\al\ : \al s \in \mathbb{N}, \al \leq \al_0 \}$ and 
$\T (s,\al) = \{ t \in \mathbb{N} \ : \ 2(s-1)+\al s\leq t \leq 
\min\{2(M-N+s),Ls\}, \ 
\mbox{$t-\al s$ is even} \}$.
We now apply (\ref{bound1}), and its corrollary, to obtain: 
\begin{eqnarray}
\tilde{\ex} [X_s] &=& O(1) {N \choose s } \sum_{\al \in \D (s ,\al_0)} 
\sum_{t \in T(s,\al)} {M \choose (t-\al s)/2}^{-1} \ {t \choose \al s} 
\nonumber \\
&=& O(1) \sum_{\al \in \D (s ,\al_0)} \sum_{t \in T(s,\al)} 
\left(\frac{te}{\al s} \right)^{\al s}M e^{-\beta s},
\nonumber
\end{eqnarray}
since $t$ is bounded by $Ls$, for small enough $\alpha$,
$({te \over \alpha s})^{\alpha s} \le e^{\beta s \over 2}$,
and each term in the above sum is at most $Me^{-\beta s \over 2}$. Further the sum has $O(s^2)$ terms. 
So, for $a_0>0$ sufficiently small and $s \geq C_1 \ln n$, where $C_1=C_1(d)$ is a 
sufficiently large constant the sum is $o({1 \over n})$. 
\end{proofof}

\section{Long Induced Paths} \label{LongPaths}

In this section we 
give the 
\vskip0.2cm

\begin{proofof} {\bf Lemma \ref{np-mix-time}} 
We focus on the  set ${\cal S}$ of those paths of $G_{n,p}$ 
all of whose internal vertices have degree two but whose 
endpoints have degree greater than 2.
Note that some of these paths are edges. 
We let ${\cal S}^\prime$ be the mulitset consisting of the interiors
of these paths. ${\cal S}^\prime$ is a multiset because it may
(indeed a.a.s does) contain many empty paths. 

The expected number of paths in ${\cal S}$ of length $i$ 
is at most $${n! \over (n-i)!}p^{i-1}(1-p)^{(i-2)(n-i)}$$
and at least $${1 \over 9} {n! \over (n-i)!}p^{i-1}(1-p)^{in}.$$
It follows that the expected number of paths of length 
$\ln~n \over 4d$ in ${\cal S}$ is $\omega(n^{3/4})$ whilst 
the expected number of paths of length  exceeding
$10 \ln~n \over d$ 
is $o(1)$.  The latter result
tells us that a.a.s every path of ${\cal S}$ has length at most
$10 \ln~n \over d$. 
The first result and a simple 
second moment calculation shows that a.a.s there will be at least 
$\sqrt{n}$ paths of length
$\ln~n \over 4d$  in ${\cal S}$. 

Now, we can construct an auxiliary multi-graph $G^\prime_{n,p}$ 
from $G_{n,p}$ by replacing each path in ${\cal S}$ by an edge
with the same endpoints (so we delete the vertices on the 
interior of these paths). 
By our bounds on the size of the 
paths in ${\cal S}$,
we know that a.a.s $G^\prime_{n,p}$ has exactly one component of
size $\Omega({n \over \ln~n})$ and this corresponds to $H_{n,p}$. 

We would like to say that given ${\cal S}$ and  $G^\prime_{n,p}$ 
we can generate $G_{n,p}$ 
by taking a uniform bijection between the paths of ${\cal S}$ 
and the edges of $G^\prime_{n,p}$ whose endpoints have degrees bigger
than 2. This may not  be true, because
$G^\prime_{n,p}$ may have multiple edges, 
and we cannot map two edges of $G^\prime_{n,p}$
to single edge paths in ${\cal S}$. 

For the moment, we assume that it is true and that
$G_{n,p}$ has $O(n\ln~n)$ 
edges,
there are at least $\sqrt{n}$ paths of length
$\ln~n \over 4d$  in ${\cal S}$, and 
the component of $G^\prime_{n,p}$ corresponding to  $H_{n,p}$ has 
$\Omega({n \over \ln~n})$ edges. Then, 
the probability that no path of ${\cal S}$ in $H_{n,p}$ 
has length $\ln~n \over 4d$  is at 
most $\left(1 - \Omega\left({1 \over \ln^2~n}\right)\right)^{\sqrt{n}}$
which is $o(1)$. 
But we have seen that these three conditions a.a.s hold, so we are done
if $G^\prime_{n,p}$  is a simple graph.

If $G^\prime_{n,p}$ is not simple, instead of considering ${\cal S}$ we consider
the subset ${\cal S}^*$ consisting of those paths of ${\cal S}$ whose 
endpoints are not joined by another path of ${\cal S}$. We let 
$G^*_{n,p}$ be the graph obtained by replacing the paths of 
$S^*$ by edges. Now, given $G^*_{n,p}$ and $S^*$, we obtain $G_{n,p}$
by taking a uniform bijection between the paths of $S^*$ and the edges of 
$G^*_{n,p}$ whose endpoints both have degree at least three and are not joined 
by a path all of whose internal vertices have degree two. An easy first 
moment calculation shows that $|S-S^*|$ is a.a.s $o(|S|)$ so we can apply the
above argument to $S^*$ to prove our result. We omit the details.
\end{proofof}

\section{Concluding Remarks} 

We have investigated some geometric properties of the giant component of supercritical
random graphs. In particular, we analysed the edge expansion of connected
subsets of various sizes, concluding that it is the sets of order
$\frac{\ln~n}{d}$ that have small expansion. As a consequense, these sets delay the mixing of
a random walk. However, as the average degree grows, these shrink and the
mixing time is determined by a global parameter which is the
diameter. This was shown in~\cite{Hild}, for average degree at least
$\ln^2~n$. In a forthcoming paper~\cite{Evol2}, we give a more detailed analysis of this
situation for degrees between $\sqrt{\ln~n \ln \ln~n}$ and $\ln^2~n$.

\end{document}